\documentclass[11pt]{amsart}

\usepackage{amsmath,amsfonts,amssymb}
\usepackage{dsfont}
\usepackage{array,delarray}
\usepackage[dvips]{graphicx}
\usepackage{xspace}
\usepackage{color}
\usepackage[square,sort&compress]{natbib}

\graphicspath{{./}{./figs/}}

\def\ds {\displaystyle}

\newcommand{\ff}{\mathbf{F}}

\newtheorem{theorem}{Theorem}[section]   
\newtheorem{corollary}[theorem]{Corollary}  
\theoremstyle{definition}
\newtheorem{definition}[theorem]{Definition}   
\theoremstyle{remark}
\newtheorem{example}[theorem]{Example}        
\numberwithin{equation}{section}     

\long\def\Comment#1{
{
 \endIgnore{\relax} 
 }}
\def\endIgnore{\relax}

\begin{document}

\title
{On the algebraic geometry of polynomial dynamical systems}

\author[A.S.~Jarrah]{Abdul~S.~Jarrah}
\address{Virginia Bioinformatics Institute\\
Virginia Polytechnic Institute and State University} 
\email[Abdul~S.~Jarrah]{ajarrah@vbi.vt.edu}

\author[R.~Laubenbacher]{Reinhard~Laubenbacher}
\address{Virginia Bioinformatics Institute\\
Virginia Polytechnic Institute and State University} 
\email[Reinhard~Laubenbacher]{reinhard@vbi.vt.edu}

\date{\today}


\begin{abstract}
This paper focuses on polynomial dynamical systems over finite
fields.  These systems appear in a variety of contexts, in
computer science, engineering, and computational biology, for
instance as models of intracellular biochemical networks.  It is
shown that several problems relating to their structure and
dynamics, as well as control theory, can be formulated and solved
in the language of algebraic geometry.

\end{abstract}

\maketitle

\section{Introduction}
The study of the dynamics of polynomial maps as time-discrete
dynamical systems has a long tradition, in particular polynomials
over the complex numbers leading to fractal structures. An example
of more recent work using techniques from algebraic geometry
includes the study of the algebraic and topological entropy of
iterates of monomial mappings $f = (f_1,\ldots ,f_n)$ over subsets
of $\mathbf C^n$ \cite{hasselblatt:07}. The study of iterates of
polynomial mappings (or, more generally, rational maps) over the
$p$-adic numbers originally arose in Diophantine geometry
\cite{Call:93}. For more recent work see, e.g.,
\cite{vivaldi:2001}.  Finally, there is a long tradition of
studying the iterates of polynomial maps $f:\mathbf
F_q\longrightarrow \mathbf F_q$ over finite fields, primarily
using techniques from combinatorics and algebraic number theory
(see, e.g., \cite{LM1, LM2}.

The general case of \emph{finite dynamical systems}
$$
f =  (f_1,\ldots ,f_n):\mathbf F_q^n\longrightarrow \mathbf F_q^n,
$$
with $f_i\in \mathbf F_q[x_1,\ldots ,x_n]$ has long been of
interest in the special case $q=2$, which includes Boolean
networks and cellular automata.  They are of considerable interest
in engineering and computer science, for instance.  Since the
1960s they are also being used increasingly as models of diverse
biological and social networks.

The underlying mathematical questions in these studies for
different fields are similar.  They relate primarily to
understanding the long-term behavior of the iterates of the
mapping.
In the case of monomial mappings, the matrix of exponents is
usually the right mathematical object to analyze by using
different methods based on the ground field. Generally, one would
like to be able to use some feature of the structure of the
coordinate functions $f_i$ to deduce properties of the structure
of the phase space of the system $f$ which represents the system's
dynamics. For finite systems, i.e., polynomial dynamical systems
$f$ over a finite field $k$, the phase space $\mathcal P(f)$ has
the form of a directed graph with vertices the elements of $k^n$.
There is an edge $\mathbf v\longrightarrow \mathbf w$ in $\mathcal
P(f)$ if $f(\mathbf v)=\mathbf w$.  Since the graph has finitely
many vertices it is easy to see that each component has the
structure of a directed \emph{limit cycle}, with trees
(\emph{transients}) feeding into each of the nodes in the cycle.

\begin{figure}[ht]
\centerline{ \raise10pt\hbox{ \framebox{
\includegraphics[width=0.1\textwidth]{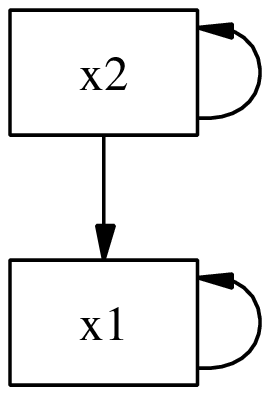}
} } \qquad \qquad \framebox{
\includegraphics[width=.4\textwidth]{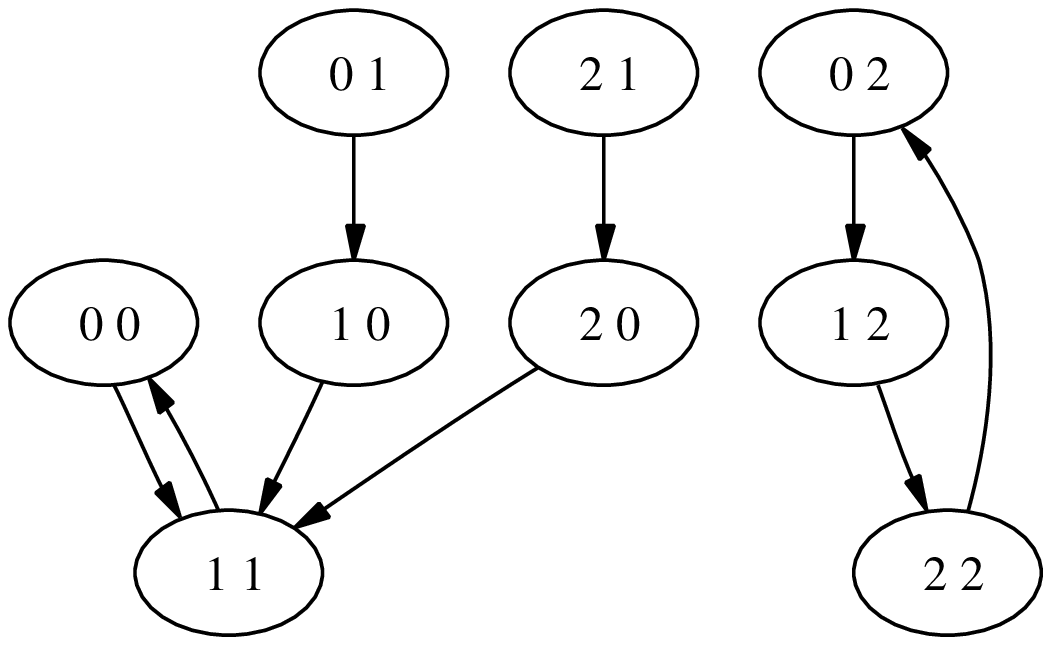}
} } \caption{The dependency graph (left) and the phase space
$\mathcal{P}(f)$ (right) of the finite dynamical system in Example
\ref{example1}.} \label{fig1}
\end{figure}

\begin{example}\label{example1}
Let $f : \ff_3^2 \longrightarrow \ff_3^2$ be given by $f(x_1,x_2)
= (1-x_1x_2, 1+2x_2)$. The phase space of $f$ has two components,
containing two limit cycles: one of length two and one of length
three. See Figure \ref{fig1} (right). The dependency relations
among the variables are encoded in the dependency graph in Figure
\ref{fig1} (left).
\end{example}

In this case the information of interest is the number of
components, the length of the limit cycles, and, possibly, the
structure of the transient trees.  It is also of interest to study
the sequence of iterates $f^r$ of $f$.

In recent years interest in polynomial dynamical systems over
general finite fields has arisen also because they are useful as
multistate generalizations of Boolean network models for
biochemical networks, such as gene regulatory networks.  In
particular, the problem of network inference from experimental
data can be formulated and solved in this framework.  The tools
used are very similar to those developed for the solution of
problems in the statistical field of design of experiments and
which led directly to the new field of algebraic statistics. (This
connection is made explicit in \cite{LS:des_exp}.) Analyzing
finite dynamical systems models of molecular networks requires
tools that provide information about network dynamics from
information about the system. Since the systems arising in
applications can be quite high-dimensional it is not feasible to
carry out such an analysis by exhaustive enumeration of the phase
space. Thus, the studies of polynomial dynamics mentioned earlier
become relevant to this type of application. Finally, an important
role of models in biology and engineering is their use for the
construction of controllers, be it for drug delivery to fight
cancer or control of an airfoil. Here too, approaches have been
developed for control systems in the context of polynomial
dynamical systems.

Here we describe this circle of ideas and some existing results.
Along the way we will mention open problems that arise. We believe
that algebraic geometry can play an essential role as a conceptual
and computational tool in this application area. Consistent with
our specific research interests we will restrict ourselves to
polynomial dynamical systems over finite fields.

\section{Finite dynamical systems}

Let $k$ be a finite field and let
$$
f = (f_1, \ldots ,f_n):k^n\longrightarrow k^n
$$
be a mapping.  It is well-known that the coordinate functions
$f_i$ can be represented as polynomials \cite[p. 369]{LN}, and
this representation is unique if we require that the degree of
every variable in every term of the $f_i$ is less than $|k|$.  If
$k$ is the field with two elements, then it is easily seen that
$f$ represents a Boolean network.  Conversely, every Boolean
network can be represented in polynomial form. Hence, polynomial
dynamical systems over a finite field, which we shall call
\emph{finite dynamical systems} include \emph{all} time-discrete
dynamical systems over finite fields. Iteration of $f$ generates
dynamics which is captured by the phase space $\mathcal P(f)$ of
$f$, defined above.  We will first focus on the problem of
inferring information about the structure of $\mathcal P(f)$ from
the polynomials $f_i$.

In principle, a lot of information about $\mathcal P(f)$ can be
gained by solving systems of polynomial equations.  For instance,
the fixed points of $f$ are precisely the points on the variety
given by the system
$$
f_1(x_1,\ldots ,x_n)-x_1=0,\ldots , f_n(x_1,\ldots ,x_n)-x_n=0.
$$
The points of period 2 are the points on the variety given by the
system $f_i^2(x_1,\ldots , x_n)-x_i=0, i=1,\ldots ,n$, and so on.
But often the most efficient way of solving these systems over a
finite field, in particular a small finite field, is by exhaustive
enumeration, which is not feasible for large $n$.  A more modest
goal would be to find out the number of periodic points of a given
period, that is, the number of points on the corresponding
variety. For finite fields this is also quite difficult without
the availability of good general tools, and much work remains to
be done in this direction.

The case of linear systems is the only case that has so far been
treated systematically, using methods from linear algebra, as one
might expect.  Let $M$ be the matrix of $f$.  Then complete
information about the number of components in $\mathcal P(f)$, the
lengths of all the limit cycles and the structure of the transient
trees can be computed from the invariant factors of $A$, together
with the orders of their irreducible factors \cite{He}.  See
\cite{JLSV} for an implementation of this algorithm in the
computer algebra system Singular \cite{Sing}.

There are results available for several other special families of
systems, using ad hoc combinatorial and graph-theoretical methods.
For instance, we have investigated the case where all the $f_i$
are monomials.  It is of interest to characterize monomial systems
all of whose periodic points are fixed points.  For the field with
two elements this characterization can be given in terms of the
\emph{dependency graph} of the system.  This graph has as vertices
the variables $x_1,\ldots ,x_n$.  There is a directed edge
$x_i\rightarrow x_j$ if $x_i$ appears in $f_j$; see Figure
\ref{fig1}(left) for an example. Given a directed graph it can be
decomposed into strongly connected components, with each component
a subgraph in which each vertex can be reached from every other
vertex by a directed path.  Associated to a strongly connected
graph we have its \emph{loop number} which is defined as the
greatest common divisor of the lengths of all directed loops in
the graph based at a fixed vertex.  (It is easy to see that this
number is independent of the vertex chosen.) The loop number is
also known as the \emph{index of imprimitivity} of the graph.

It is shown in \cite{colon:04} that the length of any cycle in the
phase space of a monomial system $f$ must divide the loop number
of its dependency graph. In particular, we have the following
result about fixed-point monomial systems.

\begin{theorem}[\cite{colon:04}]  All periodic points of a monomial
system $f$ are fixed points if and only if each strongly connected
component of the dependency graph of $f$ has loop number 1.
\end{theorem}

The study of monomial systems over general finite fields can be
reduced to studying linear and Boolean systems \cite{CJLS}.  There
are other interesting families of systems whose dynamics has been
studied. One such family is that of Boolean networks constructed
from \emph{nested canalyzing functions}. These were introduced and
studied in \cite{kpst03,kauffman:04}. The context is the use of
Boolean networks as models for gene regulatory networks initiated
by S. Kauffman \cite{Kauffman:69}. We first recall the definitions
of canalyzing and nested canalyzing functions from \cite{kpst03}.

\begin{definition}
A \emph{canalyzing function} is a Boolean function with the
property that one of its inputs alone can determine the output
value, for either ``true" or ``false" input. This input value is
referred to as the \emph{canalyzing value}, while the output value
is the \emph{canalyzed value}.
\end{definition}

\begin{example}
The function $f(x,y) = xy$ is a canalyzing function in the
variable $x$ with canalyzing value 0 and canalyzed value 0.
However, the function $f(x,y) = x+y$ is not canalyzing in either
variable.
\end{example}

Nested canalyzing functions are a natural specialization of
canalyzing functions. They arise from the question of what happens
when the function does not get the canalyzing value as input but
instead has to rely on its other inputs. Throughout this paper,
when we refer to a function of $n$ variables, we mean that $f$
depends on all $n$ variables.  That is, for $1 \leq i \leq n$,
there exists $(a_1,\dots,a_n) \in \ff_2^n$ such that
$f(a_1,\dots,a_{i-1},a_i,a_{i+1},\dots,a_n) \neq
f(a_1,\dots,a_{i-1},1+a_i,a_{i+1},\dots,a_n)$.

\begin{definition}\label{def-ncf}
A Boolean function $f$ in $n$ variables is a  \emph{nested
canalyzing function}(NCF) in the variable order $x_1, x_2, \ldots,
x_n$ with canalyzing  input values $a_1,\dots, a_n$ and canalyzed
output values $b_1,\dots,b_n$, respectively,  if it can be
expressed in the form
\begin{equation*} \label{ncf kauff}
f(x_1,x_2,\ldots, x_n) =
\begin{cases}
b_1 & ~{\rm if}~ x_1 = a_1, \\
               b_2 & ~{\rm if}~ x_1 \ne a_1 ~{\rm and}~ x_2 = a_2, \\
               b_3 & ~{\rm if}~ x_1 \ne a_1 ~{\rm and}~ x_2 \ne a_2 ~{\rm and}~ x_3 = a_3, \\
               \vdots & \hspace{1cm} \vdots \\
               b_n & ~{\rm if}~ x_1 \ne a_1 ~{\rm and}~  \cdots ~{\rm and}~ x_{n-1} \ne a_{n-1} ~{\rm and}~ x_n = a_n, \\
               b_{n}+1 & ~{\rm if}~ x_1 \ne a_1 ~{\rm and}~ \cdots ~{\rm and}~ x_n \ne a_n.
\end{cases}
\end{equation*}
\end{definition}

\begin{example}
The function $f(x,y,z) = x(y-1)z$ is nested canalyzing in the
variable order $x,y,z$ with canalyzing values 0,1,0 and canalyzed
values 0,0,0, respectively. However, the function $f(x,y,z,w) =
xy(z+w)$ is not a nested canalyzing function because if $x \ne 0$
and $y \ne 0$, then the value of the function is not constant for
any input values for either $z$ or $w$.
\end{example}

It is shown in \cite{kauffman:04} through extensive computer
simulations that Boolean networks constructed from nested
canalyzing functions show very stable dynamic behavior, with short
transient trees and a small number of components in their phase
space.  It is this type of dynamics that gene regulatory networks
are thought to exhibit.  It is thus reasonable to use this class
of functions preferentially in modeling such networks.  To do so
effectively it is necessary to have a better understanding of the
properties of this class, for instance, how many nested canalyzing
functions with a given number of variables there are.

In \cite{JLR} a parametrization of this class is given as follows.
The first step is to view Boolean functions as polynomials using
the following translation:
$$
x \wedge y = x\cdot y, \, x \vee y = x+y+xy, \, \neg x = x+1.
$$

It is shown in \cite{JLR} that the ring of Boolean functions
is isomorphic to the quotient ring $R = \ff_2[x_1,\dots,x_n]/I$,
where $I = \langle x_i^2-x_i : 1 \leq i \leq n \rangle$. Indexing
monomials by the subsets of $[n]:=\{1,\ldots ,n\}$ corresponding
to the variables appearing in the monomial, we can write the
elements of $R$ as
\begin{equation*}
R = \{\ds \sum_{S \subseteq [n]} c_S \prod_{i \in S} x_i \, : \,
c_S \in \ff_2\}.
\end{equation*}
As a vector space over $\ff_2$, $R$ is isomorphic to $\ff_2^{2^n}$
via the correspondence
\begin{equation*}\label{corresp}
R \ni \ds \sum_{S \subseteq [n]} c_S \prod_{i \in S} x_i
\longleftrightarrow (c_{\emptyset},\dots, c_{[n]}) \in
\ff_2^{2^n},
\end{equation*}
for a given fixed total ordering of all square-free monomials.
 That is, a polynomial function corresponds to the vector of
coefficients of the monomial summands.  The main result in
\cite{JLR} is the identification of the set of nested
canalyzing functions in $R$ with a subset $V^{ncf}$ of
$\ff_2^{2^n}$ by imposing relations on the coordinates of its
elements.

\begin{definition}
Let $\sigma$ be a permutation of the elements of the set $[n]$. We
define a new order relation $<_{\sigma}$ on the elements of $[n]$
as follows: $ \sigma (i) <_{\sigma} \sigma (j)$ if and only if $i
< j$.  Let $r_S^{\sigma}$ be the maximum element of a nonempty
subset $S$ of $[n]$ with respect to the order relation
$<_{\sigma}$. For any nonempty subset $S$ of $[n]$, the {\it
completion of S with respect to the permutation $\sigma$}, denoted
by $[r_S^{\sigma}]$, is the set  $[r_S^{\sigma}] = \{ \sigma(1),
\sigma(2), \ldots, \sigma(r_S) \}$.

Note that, if $\sigma$ is the identity permutation, then the
completion is $[r_S]$ := $\{1,2,\ldots,r_S\}$, where $r_S$ is the
largest element of $S$.
\end{definition}

\begin{theorem} \label{coeff rel sigma}\cite[Thm. 1]{AB2007}
Let $f \in R$ and let $\sigma$ be a permutation of the set $[n]$.
The polynomial $f$ is a nested canalyzing function in the order
$x_{\sigma (1)}, x_{\sigma (2)}, \ldots, x_{\sigma (n)}$, with
input values $a_{\sigma (i)}$ and corresponding output values
$b_{\sigma (i)}, 1 \le i \le n$, if and only if $c_{[n]} = 1 $
and, for any proper subset $S \subseteq [n]$,
\begin{equation*} \label{coeff formula sigma}
c_S = c_{[r_S^{\sigma}]} \prod_{\sigma (i) \in [r_S^{\sigma}]
\backslash S}  c_{[n] \backslash \{\sigma (i)\}}.
\end{equation*}
\end{theorem}

\begin{corollary}\label{variety-sigma}\cite[Cor. 1]{AB2007}
The set of points in $\ff_2^{2^n}$ corresponding to nested
canalyzing functions in the variable order $x_{\sigma (1)},
x_{\sigma (2)}, \ldots, x_{\sigma (n)}$, denoted by
$V_{\sigma}^{ncf}$, is defined by
\begin{equation*}
V_{\sigma}^{ncf} = \{(c_\emptyset,\dots,c_{[n]}) \in \ff_2^{2^n} :
c_{[n]}=1, \, c_S = c_{[r_S^{\sigma}]} \prod_{\sigma (i) \in
[r_S^{\sigma}] \backslash S} c_{[n] \backslash \{\sigma (i)\}}
\mbox {, for } S \subseteq [n]\}.
\end{equation*}
\end{corollary}

It was also shown in \cite{JLR} that
\begin{eqnarray*}
V^{ncf} &=& \bigcup_\sigma V_\sigma^{ncf}.
\end{eqnarray*}
Counting the points on this variety for small values of $n$
resulted in an integer sequence, which, with the help of the
On-Line Encyclopedia of Integer Sequences ({\tt
http://www.research.att.com/~njas/sequences/}) led to the
realization that the class of nested canalyzing functions is
identical to the class of unate cascade functions that has been
studied extensively in computer engineering literature. In
particular, using this equality, one obtains a recursive formula
for the number of nested canalyzing functions, see \cite[Corollary
2.11]{JLR}.

It is shown in \cite{AB2007} that the sets $V_\sigma^{ncf}$ are
 the irreducible components of $V^{ncf}$.  Precisely, it
is shown that for all permutations $\sigma$ on $[n]$, the ideal of
the variety $V_\sigma^{ncf}$, denoted by $I_\sigma :=
\mathds{I}(V_\sigma^{ncf})$, is a binomial prime ideal in the
polynomial ring $\overline{\ff}_2[\{c_S : S \subseteq [n]\}]$,
where $\overline{\ff}_2$ is the algebraic closure of $\ff_2$.

It remains to study this toric variety in more detail.  In
particular, a generalization of the concept of nested canalyzing
function to larger finite fields remains to be worked out.  Also,
the approach taken here can be applied to other classes of
functions important in network modeling, such as threshold
functions or monotone functions.

\section{Network inference}

Our motivation for much of the research described in the previous
section comes from our work on one of the central problems in
computational systems biology.  Due to the availability of
so-called "omics" data sets it is now feasible to think about
making large-scale mathematical models of molecular networks
involving many gene transcripts (genomics), proteins (proteomics),
and metabolites (metabolomics). One possible approach to this
problem is to build a phenomenological model based solely or
largely on the experimental data which can subsequently be refined
with additional biological information about the mechanisms of
interaction of the different molecular species.  That is, given a
data set, we are to infer a ``most likely" mathematical or
statistical model of the network that generated this data set. The
biggest challenge is that typically the network that generated the
data is high-dimensional (hundreds or thousands of molecular
species/variables) and the available data sets are typically very
small (tens or hundreds of data points).  Also, general properties
of such networks are not well-understood so that there are few
general selection criteria.  Therefore, it is not feasible to
apply many of the existing network inference methods. In this
context, two pieces of information about a molecular network are
of interest to a life scientist: the ``wiring diagram" of the
network indicating which variables causally affect which others,
and the long-term dynamic behavior of the network.  Some network
inference methods give only information about the wiring diagram,
others provide both.


 One possible
modeling framework is that of polynomial dynamical systems over
finite fields.  In the 1960s, S. Kauffman proposed Boolean
networks as good models for capturing key aspects of gene
regulation \cite{Kauffman:69}. In the 1990s so-called
\emph{logical models} were proposed by R. Thomas \cite{Thomas:89}
as models for biochemical network, which are multi-state
time-discrete dynamical systems, with both deterministic and
stochastic variants.  In \cite{LS} it was shown that the modeling
framework of polynomial dynamical systems over finite fields is a
good setting for the problem of network inference.  As pointed
out, these systems generalize Boolean networks and have many of
the same features as logical models. The network inference problem
can be formulated in this setting as follows.

Suppose that the biological system to be modeled contains $n$
variables, e.g., genes, and we measure $r+1$ time points $\mathbf
p_0, \ldots , \mathbf p_r$, using, e.g., gene chip technology,
each of which can be viewed as an $n$-dimensional real-valued
vector. The first step is to discretize the entries in the
$\mathbf p_i$ into a prime number of states, which are viewed as
entries in a finite field $k$.  If we choose to discretize into
two states by choosing a threshold, then we will obtain Boolean
networks as models. The discretization step is crucial in this
process as it represents the interface between the continuous and
discrete worlds.  Other network inference methods, such as most
dynamic Bayesian network methods, also have to carry out this
preprocessing step.  Unfortunately, there is very little work that
has been done on this problem.  We have developed a new
discretization method which is described in \cite{DVML}. It
compares favorably to other commonly used discretization methods,
using different network inference methods.

Given this data set, an {\it admissible model}
$$
f = (f_1, f_2, \ldots ,f_n): k^n\longrightarrow k^n
$$
consists of a dynamical system $f$ which satisfies the property
that
$$
f(\mathbf p_j) = (f_1(\mathbf p_j),\ldots , f_n(\mathbf p_j)) =
\mathbf p_{j+1}.
$$

The algorithm in \cite{LS} then proceeds to select such a model
$f$, which is the most likely one based on certain specified
criteria.  This is done by first reducing the problem to the case
of one variable, that is, to the problem of selecting the $f_i$
separately.  For this purpose, we compute the set of all functions
$f_i$ such that $f_i(\mathbf p_j) = \mathbf p_{j+1}^i$, that is,
all polynomial functions $f_i$ whose value on $\mathbf p_j$ is the
$i$th coordinate of $\mathbf p_{j+1}$.  This set can be
represented as the coset $f^0 + I$, where $f^0$ is a particular
such function and $I\subset k[x_1, \ldots ,x_n]$ is the ideal of
all polynomials that vanish on the given data set, also known as
the {\it ideal of points} of ${\mathbf p}_1, \ldots ,{\mathbf
p}_{r-1}$. The algorithm then chooses the normal form of an
interpolating polynomial $f^0$, based on a chosen term order.  One
drawback of the algorithm is that this choice of term order is
typically random and can of course significantly affect the form
of the model.

Modifications of the algorithm in \cite{LS} have been constructed.
The algorithm in \cite{JLSS} starts with only data as input and
computes all possible minimal wiring diagrams of polynomial models
that fit the given data and outputs a most likely one, based on
one of several possible model scoring methods. It does not depend
on the choice of a term order.  The algorithm is based on the
observation that if $f_i$ is the $i$-th coordinate function of a
model and $\mathbf p, \mathbf q$ are data points such that
$f_i(\mathbf p)\neq f_i(\mathbf q)$, then the function $f_i$ must
involve at least some of the variables corresponding to
coordinates in which $\mathbf p$ and $\mathbf q$ differ.  This
observation can be encoded in a monomial ideal $M$ which is
generated by all monomials of the form
$$
\prod_{\mathbf p_j\neq \mathbf
q_j}x_j
$$
for all pairs of points $\mathbf p\neq \mathbf q$ such that
$f_i(\mathbf p)\neq f_i(\mathbf q)$.  Now let $P =\langle
x_{j_1},\ldots , x_{j_t}\rangle$ be a minimal prime of $M$.  Then
it is not hard to see that the generators of $P$ induce a minimal
wiring diagram for $f_i$.  Conversely, every such minimal diagram
provides generators for a minimal prime of the ideal $M$.  This
algorithm has been implemented in Macaulay2 \cite{M2}. The
algorithm comes with a collection of probability distributions on
the set of minimal primes that can be used for model selection.

Another approach to the problem of dependency of the model
selection process in \cite{LS} on the chosen term order is taken
in \cite{DJSL}.  The algorithm there uses the Gr\"obner fan of the
ideal of points as a computational tool to find a most likely
wiring diagram.  It is clear that an upper bound for the number of
different models one can obtain from the algorithm in \cite{LS} by
varying the choice of term order is given by the number of cones
in the Gr\"obner fan. The algorithm in \cite{DJSL} uses
information about the frequency of appearance of the different
variables in models built for each cone to build a consensus
wiring diagram from this collection of possible models. It
computes the Deegan-Packel index of power \cite{DP} to rank
variables in order of significance.  This index was introduced in
\cite{Fetrow}, where it was computed using a Monte Carlo algorithm
to generate random term orders.  The use of the Gr\"obner fan
allows a systematic computation of this index.

Note that the model space $f^0+I$ contains {\it all} possible
polynomial functions that fit the given data.  In order to improve
the performance of model selection algorithms it would be very
useful to be able to select certain subspaces of functions that
have favorable properties as models of particular biological
systems, thereby reducing the model space. For instance, one might
consider imposing certain constraints on the structure of the
polynomials or on the resulting dynamics.  The desire to find such
constraints is what motivated the investigations described in the
previous section.  In order to select a class of polynomials with
prescribed dynamics one needs to be able to link polynomial
structure to dynamics in an algorithmic way.  Similarly, in order
to limit model selection to special classes of polynomials, such
as nested canalyzing functions in the Boolean case, one must be
able to identify efficiently the set of such functions from
the model space $f^0+I$.  This problem remains open.


\section{Control of finite dynamical systems}
Control of biological systems is an important aspect of
computational biology, ranging from the control of intracellular
biochemical pathways to chemotherapy drug delivery and control of
epidemiological processes.  In order to apply mathematical control
theory techniques it is necessary to work with a mathematical
model of the system for which control theoretic tools exist. There
is of course a very rich control-theoretic literature for systems
of differential equations.  However, the problem has also been
considered in the context of polynomial dynamical systems over
finite fields \cite{ML1, ML2, RS1}.  We briefly describe the
general setting.

The framework developed before needs to be slightly modified to
accommodate variables representing control inputs at each state
transitions as well as constraints on these inputs and on the set
of allowable initial conditions.

\begin{definition}
A \emph{controlled finite dynamical system} is a function
$$
F: k^n\times k^m\longrightarrow k^n,
$$
where the first set of variables $x_1,\ldots ,x_n$ represents the
state variables, and the second set $u_1,\ldots , u_m$ represents
the control variables.  Furthermore, we have a system of
polynomial equations
$$
Q_i(x_1,\ldots , x_n; u_1, \ldots ,u_m)=0, \,  i=1,\ldots ,r,
$$
which defines the variety of admissible control inputs, and
another system
$$
P_j(x_1,\ldots ,x_n)=0, \, j=1,\ldots ,s,
$$
which defines the variety of admissible initial conditions of the
system.  Finally, we have another polynomial system
$$
U_a(x_1,\ldots ,x_n)=0, \, a=1,\ldots ,t,
$$
which defines the variety of admissible final states.
\end{definition}

A typical optimal control problem is then stated as follows. Given
a controlled system $F$ and an admissible initialization
$(x_1,\ldots ,x_n)$, find a sequence of control inputs which drive
the system to an admissible final state, in such a way that a
suitably defined cost function is minimized.  There are several ad
hoc strategies of finding an optimal controller but much
theoretical work remains to be done in this context.

In \cite{JVDL} an application of this approach to a virus
competition problem was given, which we describe here in some
detail in order to illustrate the definitions. In a Petri dish the
center is infected with two different suitably chosen virus
strains.  It can be observed experimentally that, as the infection
spreads to the rest of the dish a clear pattern of segmentation
occurs between the two virus populations, rather than the expected
mixing of the two strains. Through addition of one or the other
type of virus over time the segmentation pattern can be
influenced.  For instance, it is possible to contain one virus
strain by strategically inoculating cells in the Petri dish with
the other strain.  A possible application of this observation
might be that one can contain the spread of a very harmful virus
by the strategic introduction of another, less harmful virus
strain.  In this context it would be of interest to develop
optimal strategies for introducing the second strain.

\begin{figure}[!htp]
\centering
\includegraphics[totalheight=6cm,angle=90]{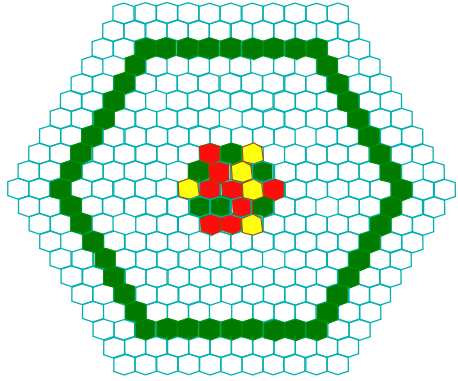}
\caption{Control objective.}  \label{control}
\end{figure}

This problem was treated in \cite{JVDL} by representing the spread
of the two virus populations as a polynomial dynamical system over
the Galois field $k=GF(4)$ as follows.  The Petri dish is
represented by 331 concentrically arranged hexagonal cells, each
of which is a variable of the system. Each cell can take on 4
possible states, corresponding to being uninfected (White),
infected by one of the two strains (Green or Red), or infected by
both (Yellow). We begin with an arbitrary initialization of the
center of the Petri dish, which is represented by the 19 innermost
cells. That is, we make an arbitrary assignment of the 4 colors to
these cells.  All remaining cells are assigned White. The goal is
to apply a series of control inputs as the infection spreads,
which prevents the Red virus from spreading to the edge of the
Petri dish.  That is, a desirable final state of the system is any
state for which the outermost ring of cells is infected only with
Green virus. See Figure \ref{control}.

The rules governing the spread of infection are as follows:
\begin{enumerate}
\item If a cell has only one infected neighbor, then it will get
the same type of infection. \item If a cell has two infected
neighbors, then one makes an assignment for the different
possibilities (for details see \cite{JVDL}).
\end{enumerate}

We use the representation $k=GF(4) := \mathbf Z_2[a] = \{0,a,
a^2,a^3\}$. The color assignment is as follows:
\[
\begin{tabular}{|l|c|}
  \hline
  Color & Field Element \\ \hline
  Green & 0 \\ \hline
  White & $a$ \\ \hline
  Yellow & $a^2$ \\ \hline
  Red & $a^3 \equiv 1$ \\ \hline
\end{tabular}
\]

If we represent the $331$ cells by the variables $x_1,\ldots
,x_{331}$, with $x_1,\ldots ,x_{19}$ representing the center cells
and $x_{271},\ldots ,x_{331}$ the cells in the outermost ring,
then the variety $U$ of admissible initializations of the model
can be described as
$$
{\bf V}(x_{20}-a,\ldots ,x_{331}-a)={\bf
V}(1-(1-(x_{20}-a)^3)\cdots (1-(x_{331}-a)^3)).
$$
(Recall that in $k=GF(4)$, we have $b^3=1$ for any nonzero $b\in
GF(4)$.)

As explained, at any point in the simulation, the cells in the
outer ring should be either green or white.  Thus, we can describe
the constraint variety $V$ as follows.  We have
\begin{eqnarray*}
{\bf V}&=&\{{\bf x}\in k^{331}| x_i(x_i-a)=0; 271\leq i\leq
331\}\\
&=&{\bf V}(1-\prod_{i=271}^{331}(1-(x_i^2-ax_i)^3)).
\end{eqnarray*}

The next step is to construct a state space model
$$
f=(f_1,\ldots ,f_{331}):k^{331}\to k^{331},
$$
of this experimental system, that is, a polynomial dynamical
system $f$ that approximates the dynamics observed in the
laboratory. The coordinate functions $f_i$ are polynomials in
$k[x_1,\ldots ,x_{331}]$ and represent the update rules of the
cells $x_i$. Since the simulated Petri dish is homogeneous, all
cells are identical, so that $f_i=f_j$ for all $i,j$.  This is
done using the algorithm in \cite{LS}.

Let $x$ represent one of the $331$ cells, and let $y_1,\ldots
,y_6$ represent its six immediate neighbors.  We compute a
symmetric polynomial $h \in k[y_1,\dots,y_6]$ which represents the
rules for the spread of the infection. That is, for a given
infection status of the neighboring six cells, the polynomial
takes on the appropriate value in $GF(4)$. At the same time, we
compute the ideal of these points $I$, that is, the set of all
polynomial functions in $k[y_1,\dots,y_6]$ that vanish at these
points. Since we are interested only in symmetric functions, let
$\bar{I}$ be the ideal of symmetric functions inside the ideal
$I$. The ideal $\bar{I}$ can easily be computed using
computational algebra methods. Thus any possible symmetric
polynomial function that can be a model of our system must be in
the set $h+\bar{I}$. We choose as a model the normal form of $h$
in the ideal $\bar{I}$, with respect to a chosen term order.  In
our case,
\begin{eqnarray*}
f &=&
\gamma_2^2+\gamma_2\gamma_1^3+a^2\gamma_1^3+a^2\gamma_1^2+a^2\gamma_1
\end{eqnarray*}
where
\begin{eqnarray*}
  \gamma_1 &=&  y_1 + \cdots + y_6, \mbox{ and }\\
  \gamma_2 &=& \sum_{i\neq j} y_iy_j.
\end{eqnarray*}
The polynomials $\gamma_1$ and $\gamma_2$ are elementary symmetric
functions in the polynomial ring $k[y_1,\dots,y_6]$. Thus, the
polynomial dynamical system $f:k^{331}\to k^{331}$ is a state
space model for our system.

The next step in formulating the optimal control problem is to
define a cost function for a given controller $g: GF(4)^{331}\to
GF(4)^{331}$.  We assume a uniform unit cost $c$ for each cell
that the controller infects with GREEN virus. Furthermore, assume
that there is an "overhead" cost $d$ attached to each
intervention, independent of the number of cells transformed. Then
the cost function for a controller $g$ is given as follows.
Suppose that $\{u_1,\ldots ,u_r\}$ is a sequence of control
inputs.  Let $c_i$ be the number of cells infected with GREEN
during control input $u_i$.  Then
$$
C(g,\{u_i\})=\sum_{i=1}^{r}(c\cdot c_i+|\delta_{c_i,0}-1|\cdot d)
.
$$
To find a controller $g$ that minimizes the cost function $C(g)$
amounts to finding a way to control the system by transforming a
minimal number of cells.

The paper \cite{JVDL} contains the construction and implementation
of a suitable controller that was experimentally verified to
accomplish the stated goal of containing one of the two strains.
The authors were not able to show that the chosen control strategy
was optimal, however.

It is well-known that optimal control of non-linear systems is
difficult and few general tools are available. A common ad hoc
approach is to ``work backwards" from a desired end state and
reconstruct control inputs in this way. One might hope that a
suitable formulation of the problem in the language of algebraic
geometry might allow the use of new tools from this field.

\section{Discussion}
It was shown in this paper that several problems about polynomial
dynamical systems over finite fields can be formulated and
possibly solved in the context of algebraic geometry and
computational algebra.  In particular, the problem of relating the
structure of the defining polynomials to the resulting dynamics
can be approached in this way.  Likewise, the inference of a
system from a given partial data set is amenable to a solution
using tools from algebraic geometry.  Finally, the beginnings of a
control theory for such systems is expressed in the language of
ideals and varieties.

However, it is clear that the results presented here barely
scratch the surface of the problems and of what could be
accomplished with the use of more sophisticated algebraic
geometric tools.


\section{Acknowledgements}
 The authors were partially supported by NSF Grant
 DMS-0511441 and NIH Grant R01 GM068947-01.

\bibliographystyle{siam}
\bibliography{IMA_Appl_Alg_Geom}

\begin{thebibliography}{10}

\bibitem{Call:93}
{\sc G.~Call and J.~Silverman}, {\em Canonical height on varieties with
  morphisms}, Compositio Math., 89 (1993), pp.~163--205.

\bibitem{CJLS}
{\sc O.~Col\'on-Reyes, A.~Jarrah, R.~Laubenbacher, and B.~Sturmfels}, {\em
  Monomial dynamical systems over finite fields}, Complex Systems, 16 (2006),
  pp.~333--342.

\bibitem{colon:04}
{\sc O.~Colon-Reyes, R.~Laubenbacher, and B.~Pareigis}, {\em Boolean monomial
  dynamical systems}, Annals of Combinatorics, 8 (2004), pp.~425--439.

\bibitem{DP}
{\sc J.~Deegan and E.~Packel}, {\em A new index for simple $n$-person games},
  Int. J. Game Theory, 7 (1978), pp.~113--123.

\bibitem{DVML}
{\sc E.~Dimitorva, P.~Vera-Licoa, J.~McGee, and R.~Laubnebahcer}, {\em
  Discretization of time series data}, 2007.
\newblock Submitted.

\bibitem{DJSL}
{\sc E.~Dimitrova, A.~Jarrah, B.~Stigler, and R.~Laubenabcher}, {\em A
  {G}roebner fan-based method for biochemical network}, in ISSAC Proceedings,
  2007.

\bibitem{Fetrow}
{\sc L.~D. S.~T. E.~Allen, J.~Fetrow and D.~John}, {\em Algebraic dependency
  models of protein signal transduction networks from time series data}, J.
  Theor. Biol., 238 (2006), pp.~317--330.

\bibitem{M2}
{\sc D.~Grayson and M.~Stillman}, {\em {\sc Macaulay 2}, a software system for
  research in algebraic geometry}.
\newblock World Wide Web.
\newblock {\texttt{http://www.math.uiuc.edu/Macaulay2}}.

\bibitem{Sing}
{\sc G.-M. Greuel, G.~Pfister, and H.~Sch{\"o}nemann}, {\em Singular 2.0}, a
  computer algebra system for polynomial computations, Centre for Computer
  Algebra, University of Kaiserslautern, 2001.
\newblock \texttt{http://www.singular.uni-kl.de}.

\bibitem{hasselblatt:07}
{\sc B.~Hasselblatt and J.~Propp}, {\em Degree growth of monomial maps}, 2006.
\newblock arXiv:Math.DS/0604521 v2.

\bibitem{He}
{\sc A.~Hern\'andez-Toledo}, {\em Linear finite dynamical systems},
  Communications in Algebra, 33 (2005), pp.~2977--2989.

\bibitem{AB2007}
{\sc A.~Jarrah and R.~Laubenbacher}, {\em Discrete models of biochemical
  networks: The toric variety of nested canalyzing functions}, in Algebraic
  Biology, H.~Anai, K.~Horimoto, and T.~Kutsia, eds., no.~4545 in LNCS,
  Springer, 2007, pp.~15--22.

\bibitem{JLSS}
{\sc A.~Jarrah, R.~Laubenbacher, B.~Stigler, and M.~Stillman}, {\em
  Reverse-engineering of polynomial dynamical systems}, Advances in Applied
  Mathematics, 39.

\bibitem{JLSV}
{\sc A.~Jarrah, R.~Laubenbacher, M.~Stillman, and P.~Vera-Licona}, {\em An
  efficient algorithm for the phase space structure of linear dynamical systems
  over finite fields}.
\newblock Submitted, 2007.

\bibitem{JLR}
{\sc A.~Jarrah, B.~Raposa, and R.~Laubenbacher}, {\em Nested canalyzing, unate
  cascade, and polynomial functions}, Physica D, 233 (2007), pp.~167--174.

\bibitem{JVDL}
{\sc A.~Jarrah, H.~Vastani, K.~Duca, and R.~Laubenbacher}, {\em An optimal
  control problem for \it{in vitro} virus competition}, in 43rd IEEE Conference
  on Decision and Control, December 2004.
\newblock Invited paper.

\bibitem{kpst03}
{\sc S.~Kauffman, C.~Peterson, B.~Samuelsson, and C.~Troein}, {\em Random
  boolean network models and the yeast transcriptional network}, Proc. Natl.
  Acad. Sci. USA., 100 (2003), pp.~14796--9.

\bibitem{kauffman:04}
{\sc S.~Kauffman, C.~Peterson, B.~Samuelsson, and C.~Troein}, {\em {Genetic
  networks with canalyzing Boolean rules are always stable}}, PNAS, 101 (2004),
  pp.~17102--17107.

\bibitem{Kauffman:69}
{\sc S.~A. Kauffman}, {\em Metabolic stability and epigenesis in randomly
  constructed genetic nets}, Journal of Theoretical Biology, 22 (1969),
  pp.~437--467.

\bibitem{LS}
{\sc R.~Laubenbacher and B.~Stigler}, {\em A computational algebra approach to
  the reverse-engineering of gene regulatory networks}, J Theor Bio, 229
  (2004), pp.~523--537.

\bibitem{LS:des_exp}
\leavevmode\vrule height 2pt depth -1.6pt width 23pt, {\em Design of
  experiments and biochemical network inference}, in Algebraic and Geometric
  Methods in Statistics, R.~E. Gibilisco~P., ed., Cambridge University Press,
  Cambridge, 2007.

\bibitem{LM1}
{\sc L.~Lidl and G.~Mullen}, {\em When does a polynomial over a finite field
  permute the elements of the field?}, American Mathematical Monthly, 95
  (1988), pp.~243--246.

\bibitem{LM2}
\leavevmode\vrule height 2pt depth -1.6pt width 23pt, {\em When does a
  polynomial over a finite field permute the elements of the field?,{II}},
  American Mathematical Monthly, 100 (1993), pp.~71--74.

\bibitem{LN}
{\sc R.~Lidl and H.~Niederreiter}, {\em Finite Fields}, Cambridge University
  Press, New York, 1997.

\bibitem{ML1}
{\sc H.~Marchand and M.~LeBorgne}, {\em On the optimal control of polynomial
  dynamical systems over $\zz/p\zz$}, in Fourth Workshop on Discrete Event
  Systems, IEEE, Cagliari, Italy, 1998.

\bibitem{ML2}
\leavevmode\vrule height 2pt depth -1.6pt width 23pt, {\em Partial order
  control of discrete event systems modeled as polynomial dynamical systems},
  in IEEE International conference on control applications, Trieste, Italy,
  1998.

\bibitem{vivaldi:2001}
{\sc J.~Pettigrew, J.~Roberts, and F.~Vivaldi}, {\em Complexity of regular
  invertible $p$-adic motions}, Chaos, 11 (2001), pp.~849--857.

\bibitem{RS1}
{\sc L.~Reger and K.~Schmidt}, {\em Aspects on analysis and synthesis of linear
  discrete systems over the finite field $gf(q)$}, in Proc. European Control
  Conference ECC2003, Cambridge University Press, 2003.

\bibitem{Thomas:89}
{\sc R.~Thomas and R.~D'Ari}, {\em Biological Feedback}, CRC Press, 1989.

\end{thebibliography}

\end{document}